\documentclass[12pt]{article}
\usepackage{amsmath,amssymb,amsthm, amscd}

\begin{document}

\begin{center} 
{\bf Notes on Poisson deformations of symplectic varieties},  
\vspace{0.2cm}

- Corrections and Addendum to the article ``Flops and Poisson deformations of symplectic varieties, Publ. Res. Inst. Math. Sci. {\bf 44} (2008) 259 - 314'' - 
\vspace{0.2cm}

{\bf Yoshinori Namikawa}

\end{center}
\vspace{0.2cm}

\S 1. {\bf Introduction}.  

We have proved some basic properties of the Poisson deformation of a symplectic variety $(X, \omega)$ in [Na]. 
We call an affine symplectic variety $(X, \omega)$ {\em conical} 
if $X$ admits a good $\mathbf{C}^*$-action and $\omega$ is homogeneous with weight $l > 0$.  
Assume that a symplectic variety $(X, \omega)$ is a crepant partial projective resolution of a 
conical symplectic variety $(Y, \omega_Y)$ and $X$ admits only terminal singularities. Then the $\mathbf{C}^*$-action on $Y$ extends to 
a $\mathbf{C}^*$-action.on $X$.  
Lemma 20 of [Na] claims that this $\mathbf{C}^*$-action indces a natural $\bf{C}^*$-action on the universal formal Poisson deformation 
of $X$. But the proof given there was insufficient. One purpose of this note is to 
give a precise argument to this part. 

Another purpose of this note is to give an alternative proof of Lemma A.8 of [Na], where the original proof contains an error.  

In the following we freely use the terminology of [Na] and [Na 2]. .   
\vspace{0.2cm}

\S 2. {\bf The universal Poisson deformations and $\mathbf{C}^*$-actions}  
\vspace{0.2cm} 

Let $(X, \omega)$ be a symplectic variety. Assme that one of the following holds: 
\vspace{0.2cm}

(i) $(X, \omega)$ is a conical symplectic variety with $wt(\omega) = l > 0$. 
\vspace{0.2cm} 

(ii) $(X, \omega)$ is a crepant partial projective resolution of a 
conical symplectic variety $(Y, \omega_Y)$ with $wt(\omega_Y) = l > 0$ and $X$ admits only terminal singularities.
\vspace{0.2cm}

Let  $(X, \{\:, \:\})$ be the Poisson structure determined by $\omega$.  
Then the Poisson deformation functor $\mathrm{PD}_X$ is prorepresentable by Corollary 2.5 of [Na 2].  Let $R$ be the prorepresentable hull of $\mathrm{PD}_X$. 
Putting $R_n := R/m_R^{n+1}$ and $S_n := \mathrm{Spec}R_n$, let $\{(X_n, \{\:, \}_n)\}$ be 
the universal formal Poisson deformation. Here $(X_n, \{\:, \}_n)$ is a  Poisson deformation of $(X, \{\:, \:\})$ over $S_n$, and  
$(X_{n+1}, \{\:, \:\}_{n+1})$ is an extension of $(X_n, \{\:, \}_n)$ to $S_{n+1}$ for each $n$.  

\begin{equation} 
\begin{CD}
X_0 := X @>>> X_1 @>>> ... @>>> X_n  @>>> ... \\ 
@VVV @VVV @VVV @VVV @VVV\\ 
S_0 @>>> S_1 
@>>> ... @>>> S_n @>>> ...         
\end{CD} 
\end{equation}  
\vspace{0.2cm}

We shall prove the following. The proof is a modification of the argument by Rim [R]. 
\vspace{0.2cm}

{\bf Theorem (2.1)}.  {\em Assume that $R$ prorepresents $\mathrm{PD}_{(X, \{\:, \:\})}$ and 
$\{(X_n, \{\:, \}_n)\}$ is the universal formal  Poisson deformation of $(X, \{\:, \:\})$. 
Then $\{X_n\} \to \{S_n\}$ admits a natural $\mathbf{C}^*$-action induced by the $\mathbf{C}^*$-action on 
$X$. Namely, we have $\mathbf{C}^*$-actions on $X_n$ and $S_n$ with the following properties:   

(i) The $\mathbf{C}^*$-action on $X_0$ coincides with the original one, and the commutative diagram above is $\mathbf{C}^*$-equivariant. 

(ii) Denote by $(\phi_{\sigma})_n$ the automorphism of $X_n$ determined by $\sigma \in \mathbf{C}^*$. Then 
$(\phi_{\sigma})_n: (X_n, \sigma^{-l}\{\:, \:\}_n) \to (X_n, \{\:, \:\}_n)$ is a Poisson isomorphism. }
\vspace{0.2cm}

{\bf Proof}.  

(a) The ${\bf C}^*$-action on the base space: 

Let $\sigma: X_0 \to X_0$ be the 
action of $\sigma \in \mathbf{C}^*$ on $X_0$. 
Then, for the Poisson deformation $i: (X_0, \{\:, \:\}_0) \to (X_n, \{\:, \:\}_n)$ over $S_n$, we have a new Poisson deformation $i\circ
\sigma^{-1}: (X_0, \{\:, \:\}_0) \to (X_n, \sigma^{-l}\{\:, \:\}_n)$ over 
$S_n$. By the semiuniversality of $\{(X_n, \{\:, \:\}_n)\}$, 
there exists a morphism $\sigma_n: S_n \to S_n$ and, the Poisson deformation $i\circ\sigma^{-1}: (X_0, \{\:, \:\}_0) \to (X_n, \sigma^{-l}\{\:, \:\}_n)$ is equivalent to the Poisson deformation obtained by pulling back $i: (X_0, \{\:, \:\}_0) \to (X_n, \{\:, \:\}_n)$ by $\sigma_n$. 
In other words, there exists a Poisson isomorphism $(\phi_{\sigma})_n$ which commutes the following diagram: 

\begin{equation} 
\begin{CD}
(X_0, \sigma^{-l} \{\:, \:\}_0) @>{\sigma}>> (X_0, \{\:, \:\}_0) \\
@V{i}VV @V{i}VV \\ 
(X_n,  \sigma^{-l}\{\:, \:\}_n) @>{(\phi_{\sigma})_n}>> (X_n, \{\:, \:\}_n) \\ 
@VVV @VVV \\ 
S_n @>{\sigma_n}>> S_n          
\end{CD} 
\end{equation}\vspace{0.1cm}
  
Since $\{(X_n, \{\:, \:\}_n)\}$ is universal, such a morphism $\sigma_n$ is detemined uniquely by $\sigma$. However, we only have the existence of the Poisson automorphism $(\phi_{\sigma})_n$, but it is not unique. Here we shall prove that 
$(\sigma \circ \tau)_n = \sigma_n \circ \tau_n$ for 
$\tau \in \mathbf{C}^*$. For this purpose, it is enough to show that the Poisson deformation $i\circ\tau^{-1}\circ\sigma^{-1}: (X_0, \{\:, \:\}) \to (X_n, (\sigma\tau)^{-l}\{\:, \:\}_n)$ is equivalent to the Poisson deformation obtained by pulling back $i: (X_0, \{\:, \:\}_0) \to (X_n, \{\:, \:\}_n)$ by $\sigma_n \circ \tau_n: S_n \to S_n$. In fact, if so, then since $(\sigma \circ \tau)_n$ has this property by definition, 
it follows that $(\sigma \circ \tau)_n = \sigma_n \circ \tau_n$ by the uniqueness.

Now, for $\tau$, there exist a Poisson isomorphism $(\phi_{\tau})_n: (X_n, \tau^{-l}\{\:, \:\}_n) \to (X_n, \{\:, \:\}_n)$ and a similar commutative daigram. Here, noting that $(\phi_{\tau})_n$ can be also regarded as the Poisson isomorphism 
$(X_n, \tau^{-l}\sigma^{-l}\{\:, \:\}_n) \to (X_n, \sigma^{-l}\{\:, \:\}_n)$, we have a commutatrive diagram  
\begin{equation} 
\begin{CD}
(X_0, \tau^{-l}\sigma^{-l}\{\:, \:\}_0) @>{\tau}>> (X_0, \sigma^{-l}\{\:, \:\}_0) \\
@V{i}VV @V{i}VV \\ 
(X_n, \tau^{-1}\sigma^{-l}\{\:, \:\}_n) @>{(\phi_{\tau})_n}>> (X_n, \sigma^{-l}\{\:, \:\}_n) \\ 
@VVV @VVV \\ 
S_n @>{\tau_n}>> S_n          
\end{CD} 
\end{equation} Composing this commutative diagram and the previous commutative diagram, 
we have a desired commuative diagram: 

\begin{equation} 
\begin{CD}
(X_0, \tau^{-l}\sigma^{-l}\{\:, \:\}_0) @>{\sigma\circ\tau}>> (X_0, \{\:, \:\}_0) \\
@V{i}VV @V{i}VV \\ 
(X_n, \tau^{-l} \sigma^{-l}\{\:, \:\}_n) @>{(\phi_{\sigma})_n\circ(\phi_{\tau})_n}>> (X_n, \{\:, \:\}_n) \\ 
@VVV @VVV \\ 
S_n @>{\sigma_n\circ\tau_n}>> S_n          
\end{CD} 
\end{equation}.\vspace{0.1cm}\\*  

Since it is obvious that $\sigma_n \vert_{S_{n-1}} = \sigma_{n-1}$ for each $n$, the diagram 
$S_0 \to S_1 \to ... \to S_n \to ...$ of the base spaces admits (compatible) $\mathbf{C}^*$-actions.

(b) The $\mathbf{C}^*$-action on the universal family $\{X_n\}$. 

As noticed above, $(\phi_{\sigma})$ is not determined unquely for $\sigma$. 
The main point is that we need to take suitable $(\phi_{\sigma})$ so that they satisfy 
$(\phi_{\sigma})_n \circ 
(\phi_{\tau})_n = (\phi_{\sigma \circ \tau})_n$. 
We shall prove this by induction on $n$. 

Now suppose that we are given a $\mathbf{C}^*$-action on $X_{n-1}$, compatible with the $\mathbf{C}^*$-action on $S_{n-1}$, and 
it satisfies condition (ii) of Theorem (2.1). We shall extends this to a $\mathbf{C}^*$-action on $X_n$.  

First, let $\mathbf{C}[t, 1/t]$ be the coordinate ring of $\mathbf{C}^*$, and introduce a Poisson structure on the 
coordinate ring $\mathbf{C}[t, 1/t] \otimes_{\mathbf C} \mathcal{O}_{X_n}$ of 
$\mathbf{C}^* \times X_n$ by $t^{-l} \otimes \{\:, \:\}_n$. Now $\mathbf{C}^* \times X_{n}$ is a Poisson scheme over $\mathbf{C}^* \times S_n$. Quite similarly, $\mathbf{C}^* \times X_{n-1}$ is a Poisson scheme over $\mathbf{C}^* \times S_{n-1}$. Suppose that the following commutative dagram is given: 
\begin{equation} 
\begin{CD}
\mathbf{C}^* \times X_{n-1} @>{\phi_{n-1}}>> X_{n-1} \\
@VVV @VVV \\ 
\mathbf{C}^* \times S_{n-1} @>>> S_{n-1}.          
\end{CD} 
\end{equation}
Here $\mathbf{C}^* \times S_{n-1} \to S_{n-1}$ is the $\mathbf{C}^*$-action on $S_{n-1}$.  

{\bf Lemma (2.2)}. {\em Assume that $\phi_{n-1}$ satisfies the following two conditions: . 

(i) $\phi_{n-1}$ determines a $\mathbf{C}^*$-action on $X_{n-1}$. 

(ii) $\phi_{n-1}: (\mathbf{C}^* \times X_{n-1}, t^{-l} \otimes \{\:, \:\}_{n-1}) 
\to (X_{n-1}, \{\:, \:\}_{n-1})$ is a Poisson morphism.  

Then $\phi_{n-1}$ extends to a Poisson morphism $\phi_n: (\mathbf{C}^* \times X_n, t^{-l} \otimes \{\:, \:\}_n)) \to (X_n, \{\:, \:\}_n)$ so that the following diagram commutes: 
\begin{equation} 
\begin{CD}
\mathbf{C}^* \times X_n @>{\phi_n}>> X_n \\
@VVV @VVV \\ 
\mathbf{C}^* \times S_n @>>> S_n           
\end{CD} 
\end{equation}
Here $\mathbf{C}^* \times S_n \to S_n$ is the $\mathbf{C}^*$-action on $S_n$.}  
\vspace{0.2cm}

The following is the meaning of this lemma. For each $\sigma \in \mathbf{C}^*$, $\phi_{n-1}$ determines 
a Poisson isomorphism $(\phi_{\sigma})_{n-1}: (X_{n-1}, \sigma^{-l}\{\:, \:\}_{n-1}) \to (X_{n-1}, \{\:, \:\}_{n-1})$. We 
can extends these to Poisson isomorphisms $(\phi_{\sigma})_n: (X_n, \sigma^{-l}\{\:, \:\}_n) \to (X_n, \{\:, \:\}_n)$.  
Further, the lemma claims that $\{(\phi_{\sigma})_n\}_{\sigma \in \mathbf{C}^*}$ is an {\bf algebraic} family. 
By assumption $\{(\phi_{\sigma})_{n-1}\}_{\sigma \in \mathbf{C}^*}$ defines a $\mathbf{C}^*$-action on $X_{n-1}$. However, $\{(\phi_{\sigma})_n\}_{\sigma \in \mathbf{C}^*}$ is not necessarily a $\mathbf{C}^*$-action on $X_n$. 
\vspace{0.1cm}

{\bf Proof of Lemma (2.2)}: To prove the lemma, it is enough to show it by replacing $X_n$ and $X_{n-1}$ respectively by $(X_n)_{reg}$ and 
$(X_{n-1})_{reg}$. In fact, let $\iota_n :(X_n)_{reg} \to X_n$
be the natural open immersion, then it holds that $\mathcal{O}_{X_n} = (\iota_n)_*\mathcal{O}_{(X_n)_{reg}}$, 
$\mathcal{O}_{\mathbf{C}^* \times X_n} = (id \times \iota_n)_*\mathcal{O}_{\mathbf{C}^* 
\times (X_n)_{reg}}$. Therefore, the Poisson morphism  
$\mathbf{C}^* \times (X_n)_{reg} \to (X_n)_{reg}$ uniquely extends to a Poisson morphism 
$\mathbf{C}^* \times X_n \to X_n$. 

By Sumihiro's theorem (cf. [KKMS], ch. I, \S 2), $X_{\mathrm{reg}}$ has a $\mathbf{C}^*$-equivariant open covering $\{U^0_i\}_{i \in I}$.  
Then $\{U_i := X_n\vert_{U^0_i}\}_{i \in I}$ is an affine open covering of 
$(X_n)_{reg}$. Put $\bar{U}_i := X_{n-1}\vert_{U^0_i}$.  By definition. $\{\bar{U}_i\}_{i \in I}$ is a 
$\mathbf{C}^*$-equivariant affine open covering of 
$(X_{n-1})_{reg}$. By restricting $\phi_{n-1}$ to $\mathbf{C}^* \times \bar{U}_i$, we get a commutative diagram 
\begin{equation} 
\begin{CD}
\mathbf{C}^* \times \bar{U}_i @>{\phi_{n-1,i}}>> \bar{U}_i \\
@VVV @VVV \\ 
\mathbf{C}^* \times S_{n-1} @>>> S_{n-1}.          
\end{CD} 
\end{equation}
Since $U_i \to S_n$ is a smooth morphism, this commutative diagram extends to a commutative diagram 
\begin{equation} 
\begin{CD}
\mathbf{C}^* \times U_i @>{\phi_{n,i}}>> U_i \\
@VVV @VVV \\ 
\mathbf{C}^* \times S_n @>>> S_n.          
\end{CD} 
\end{equation}
We want to retake each $\phi_{n,i}$ suitably, and glue them to get a global commutative diagram 
\begin{equation} 
\begin{CD}
\mathbf{C}^* \times (X_n)_{reg} @>{\phi_n}>> (X_n)_{reg} \\
@VVV @VVV \\ 
\mathbf{C}^* \times S_n @>>> S_n.           
\end{CD} 
\end{equation}
The obstruction to getting such a diagram is a 1-Cech cocycle with coeffecients
 
$\phi_{n-1}^*\Theta_{(X_{n-1})_{reg}/S_{n-1}}\otimes_{R_{n-1}}(m_R^n/m_R^{n+1})$. 
We write it $$\{\alpha_{i,j}\} \in 
\prod_{i,j \in I}\Gamma(\mathbf{C}^* \times (U^0_i \cap U^0_j),  \Theta_{\mathbf{C}^* \times X_{\mathrm{reg}}/\mathbf{C}^*}\otimes_{\mathbf C}(m_R^n/m_R^{n+1})).$$  Here we notice that 
$$\phi_{n-1}^*\Theta_{(X_{n-1})_{reg}/S_{n-1}}\otimes_{R_{n-1}}(m_R^n/m_R^{n+1}) = 
\Theta_{\mathbf{C}^* \times X_{\mathrm{reg}}/\mathbf{C}^*}\otimes_{\mathbf C}(m_R^n/m_R^{n+1}).$$ 

On the other hand, the Poisson structure $\{\:, \:\}_n$ on $X_n$ determines a 2-vector $\theta_i \in 
\Gamma (U_i, \Theta_{U_i/S_n})$ on each $U_i$. Simillarly the Poisson structure on $\mathbf{C}^* \times X_n$ determines a 
2-vector $\Theta_i \in \Gamma(\mathbf{C}^* \times 
U_i, \Theta_{\mathbf{C}^* \times U_i/\mathbf{C}^* \times S_n})$ on $\mathbf{C}^* \times U_i$. 

Then the two 2-vectors $\phi_{n,i}^*\theta_i$ and $\Theta_i$ coincide if we restrict them to $\mathbf{C}^* \times \bar{U}_i$. From this it follows that 
$$\beta_i := \phi_{n,i}^*\theta_i - \Theta_i \in \Gamma(\mathbf{C}^* \times U^0_i, \wedge^2\Theta_{\mathbf{C}^* \times X_{\mathrm{reg}}/\mathbf{C}^*}\otimes_{\mathbf C}(m_R^n/m_R^{n+1})).$$    

Here let us consider the hypercohomology 
$\mathbf{H}^p(\mathbf{C}^* \times X_{\mathrm{reg}}, \Theta^{\geq 1}_{\mathbf{C}^* \times X_{\mathrm{reg}}/\mathbf{C}^*})$
of the truncated Lichnerowicz-Poisson complex 
$(\Theta^{\geq 1}_{\mathbf{C}^* \times X_{\mathrm{reg}}/\mathbf{C}^*}, \delta)$. 
Put $\mathbf{C}^* \times \mathcal{U} := \{\mathbf{C}^* \times U^0_i\}_{i \in I}$. 
Recall that it can be computed by means of the total complex of the Cech double complex 
$(C^{\cdot}(\mathbf{C}^* \times \mathcal{U}, \Theta^{\cdot}_{\mathbf{C}^* \times X_{\mathrm{reg}}/\mathbf{C}^*}); \delta, \delta_{cech})$  
\begin{equation} 
\begin{CD}
@A{\delta}AA @A{\delta}AA \\
C^0(\mathbf{C}^* \times \mathcal{U}, \Theta^2_{\mathbf{C}^* \times X_{\mathrm{reg}}/\mathbf{C}^*})   @>{\delta_{cech}}>> 
C^1(\mathbf{C}^* \times \mathcal{U}, \Theta^2_{\mathbf{C}^* \times X_{\mathrm{reg}}/\mathbf{C}^*})
 @>{-\delta_{cech}}>> \\ 
@A{\delta}AA @A{\delta}AA \\ 
C^0(\mathbf{C}^* \times \mathcal{U}, \Theta_{\mathbf{C}^* \times X_{\mathrm{reg}}/\mathbf{C}^*})   @>{-\delta_{cech}}>> 
C^1(\mathbf{C}^* \times \mathcal{U}, \Theta_{\mathbf{C}^* \times X_{\mathrm{reg}}/\mathbf{C}^*})   @>{\delta_{cech}}>>      
\end{CD} 
\end{equation} 
Here the left-most term on the bottom row has degree $1$.  

The $(\{\alpha_{i,j}\}, \{\beta_i\})$ constructed above, can be regarded as a 2-cocycle of 
$$(C^{\cdot}(\mathbf{C}^* \times \mathcal{U}, \Theta^{\cdot}_{\mathbf{C}^* \times X_{\mathrm{reg}}/\mathbf{C}^*}) ; \delta, \delta_{cech}) \otimes_{\mathbf C}(m_R^n/m_R^{n+1}).$$ 
Now we put $ob := [(\{\alpha_{i,j}\}, \{\beta_i\})] \in 
\mathbf{H}^2(\mathbf{C}^* \times X_{\mathrm{reg}}, \Theta^{\geq 1}_{\mathbf{C}^* \times X_{\mathrm{reg}}/\mathbf{C}^*}) \otimes_{\mathbf C}(m_R^n/m_R^{n+1})$. 
Since the claim of the lemma is equivalent to $ob = 0$, we need to show that $ob = 0$. 

Note that the Lichnerowicz-Poisson complex $(\Theta^{\cdot}_{\mathbf{C}^* \times X_{\mathrm{reg}}/\mathbf{C}^*}, \delta)$ is isomorphic to the de-Rham complex $(\Omega^{\cdot}_{\mathbf{C}^* \times X_{\mathrm{reg}}/\mathbf{C}^*}, d)$ and 
$(\Omega^{\cdot}_{\mathbf{C}^* \times X_{\mathrm{reg}}/\mathbf{C}^*}, d) \cong 
(\Omega^{\cdot}_{X_{\mathrm{reg}}}, d) \otimes_{\mathbf C}\mathcal{O}_{\mathbf{C}^*} \cong 
(\Theta^{\cdot}_{X_{\mathrm{reg}}}, \delta)\otimes_{\mathbf C}\mathcal{O}_{\mathbf{C}^*}$. We then have 
$$\mathbf{H}^2(\mathbf{C}^* \times X_{\mathrm{reg}}, \Theta^{\geq 1}_{\mathbf{C}^* \times X_{\mathrm{reg}}/\mathbf{C}^*}) \cong \mathbf{H}^2(X_{\mathrm{reg}}, \Theta^{\geq 1}_{X_{\mathrm{reg}}}) \otimes_{\mathbf C}\Gamma (\mathbf{C}^*, \mathcal{O}_{\mathbf{C}^*}).$$
Define the evaluation map $ev_{\sigma}: \Gamma (\mathbf{C}^*, \mathcal{O}_{\mathbf{C}^*}) \to \mathbf{C}$ at $\sigma \in \mathbf{C}^*$ by  
$ev_{\sigma}(f) = f(\sigma)$. 

We denote by $ob(\sigma)$ the image of $ob$ by the composite $$\mathbf{H}^2(\mathbf{C}^* \times X_{\mathrm{reg}}, \Theta^{\geq 1}_{\mathbf{C}^* \times X_{\mathrm{reg}}/\mathbf{C}^*}) \otimes_{\mathbf C} (m_R^n/m_R^{n+1})$$ $$\cong \mathbf{H}^2(X_{\mathrm{reg}}, \Theta^{\geq 1}_{X_{\mathrm{reg}}}) \otimes_{\mathbf C}\Gamma (\mathbf{C}^*, \mathcal{O}_{\mathbf{C}^*}) \otimes_{\mathbf C} (m_R^n/m_R^{n+1})$$ $$\stackrel{id \otimes ev_{\sigma} \otimes id}\to \mathbf{H}^2(X_{\mathrm{reg}}, \Theta^{\geq 1}_{X_{\mathrm{reg}}}) \otimes_{\mathbf C} (m_R^n/m_R^{n+1}).$$ In order to show that $ob = 0$, it is enough to show that $ob(\sigma) = 0$ 
for $\sigma$. 

By the construction of $ob$, we see that $ob(\sigma)$ is the obstruction to extending 
\begin{equation} 
\begin{CD}
(X_{n-1}, \sigma^{-l}\{\:, \:\}_n) @>{(\phi_{n-1})_{\sigma}}>> (X_{n-1}, \{\:, \:\}_{n-1}) \\
@VVV @VVV \\ 
S_{n-1} @>{\sigma_{n-1}}>> S_{n-1}          
\end{CD} 
\end{equation}
to a commutative diagram 
\begin{equation} 
\begin{CD}
(X_n, \sigma^{-l}\{\:, \:\}_n) @>>> (X_{n-1}, \{\:, \:\}_n) \\
@VVV @VVV \\ 
S_n @>{\sigma_n}>> S_n.          
\end{CD} 
\end{equation}\vspace{0.1cm}\\*
Therefore, if an extension exists, then we see that $ob(\sigma) = 0$. 
In fact, one can construct an extension of $(\phi_{n-1})_{\sigma}$ as follows. 

First, since $\{X_n, \{\:, \:\}_n\}$ is universal, we know the existence of the commtative diagram 
\begin{equation} 
\begin{CD}
(X_n, \sigma^{-l}\{\:, \:\}_n) @>{(\phi'_{\sigma})_n}>> (X_{n-1}, \{\:, \:\}_n) \\
@VVV @VVV \\ 
S_n @>{\sigma_n}>> S_n.          
\end{CD} 
\end{equation}
But, in general, it is not clear that this commutative diagram is an extension of the originally given commutative diagram. 
We only know that  $(\phi'_{\sigma})_n$ is an extension of $\sigma: (X_0, \sigma^{-l}\{\:, \:\}_0) \to 
(X_0, \{\:, \:\}_0)$. Here $(\phi_{\sigma}) \circ 
((\phi'_{\sigma})_n\vert_{X_{n-1}})^{-1}$ is an $S_{n-1}$-Poisson automorphism of $(X_{n-1}, \{\:, \:\}_{n-1})$ which restricts to the 
identity map of $X_0$. Now, since the Poisson deformation functor $\mathrm{PD}_{(X, \{\:, \:\})}$
is prorepresentable,  $(\phi_{\sigma}) \circ 
((\phi'_{\sigma})_n\vert_{X_{n-1}})^{-1}$ always extends to an $S_n$-Poisson automorphism $\psi$ of $(X_n, \{\:, \:\}_n)$. 
Then, putting $(\phi_{\sigma})_n := 
\psi \circ (\phi'_{\sigma})_n$, we see that $(\phi_{\sigma})_n$ is a desired extension.  
$\square$
\vspace{0.2cm}

Let us return to the proof of Theorem (2.1).  The problem is to find $\phi_n$ in Lemma (2.2) so that it gives a $\mathbf{C}^*$-action on 
$X_n$. To make the notation simple, we write $\eta$ for the pair of 
$\{X_0 \stackrel{i}\to X_n \stackrel{j}\to S_n\}$ and the Poisson structure $\{\:, \:\}_n$ on $X_n$. 
As already explained, for $\sigma \in \mathbf{C}^*$, we can consider the pair of $\{X_0 \stackrel{i\circ \sigma^{-1}}\to X_n \stackrel{\sigma \circ j}\to S_n\}$ and the Poisson structure $\sigma^{-l}\{\:, \:\}_n$ on $X_n$. We denote it by $\sigma\eta$. We also write simply 
$\phi_{\sigma}$ for the Poisson isomorphism $(\phi_{\sigma})_n: X_n \to X_n$ given in Lemma (2.2). 
Then $\phi_{\sigma}$ gives an $S_n$-Poisson isomorphism from $\sigma\eta$ to $\eta$, which we denote by 
$\sigma\eta \stackrel{\phi_{\sigma}}\to \eta$. 
On the other hand, $\phi_{\sigma}$ induces an $S_n$-Poisson isomorphism $\tau\sigma\eta \to \tau\eta$, which we denote by $\tau\phi_{\sigma}$. 
If, for all $\sigma, \tau \in \mathbf{C}^*$, $\phi_{\sigma}\circ \sigma\phi_{\tau} = \phi_{\sigma\tau}$, then $\{\phi_{\sigma}\}$ gives a 
$\mathbf{C}^*$-action on $X_n$. 

We write $\bar{\eta}$ for the restriction over $S_{n-1}$ of the Poisson deformation $\eta$ of  $(X_0, \{\:, \:\}_0)$ over $S_n$. 
For $\bar{\eta}$, we similarly define $\bar{\phi}_{\sigma}: \sigma\bar{\eta} \to \bar{\eta}$ and $\tau\bar{\phi}_{\sigma}: \tau\sigma\bar{\eta} \to \tau\bar{\eta}$. 
The restriction of $\phi_{\sigma}$ over $S_{n-1}$ is $\bar{\phi}_{\sigma}$, and, by assumption, we already have $\bar{\phi}_{\sigma}\circ \sigma\bar{\phi}_{\tau} = \bar{\phi}_{\sigma\tau}$. 
  
Here let $\mathrm{PAut}(\eta; id\vert_{\bar{\eta}})$ be the $\mathbf{C}$-vector space consisting of $S_n$-Poisson automorphisms of 
$\eta$ which restrict to the identity map $id$ on $\bar{\eta}$. 
For $\sigma \in \mathbf{C}^*$, $u \in \mathrm{PAut}(\eta; id\vert_{\bar{\eta}})$, we define $$^{\sigma}u := \phi_{\sigma} \circ \sigma u \circ \phi_{\sigma}^{-1}$$
$$ \eta \stackrel{\phi_{\sigma}^{-1}}\to \sigma\eta \stackrel{\sigma u}\to 
\sigma\eta \stackrel{\phi_{\sigma}}\to \eta. $$  
Then, since $\sigma\bar{u} = id$, we have $^{\sigma}u \in \mathrm{PAut}(\eta; id\vert_{\bar{\eta}})$.  Furthermore, we have 
$$^{\tau}(^{\sigma}u) = \phi_{\tau} \circ \tau(^{\sigma}u) \circ \phi_{\tau}^{-1} = 
\phi_{\tau} \circ \tau(\phi_{\sigma} \circ \sigma u \circ \phi_{\sigma}^{-1}) \circ 
\phi_{\tau}^{-1}$$
$$ =  \phi_{\tau} \circ \tau\phi_{\sigma} \circ (\tau\sigma)u \circ \tau\phi^{-1}_{\sigma} 
\circ \phi^{-1}_{\tau} = \phi_{\tau} \circ \tau\phi_{\sigma} \circ (\tau\sigma)u \circ 
(\tau \phi_{\sigma})^{-1} \circ \phi^{-1}_{\tau} $$
$$ = \phi_{\tau\sigma} \circ (\tau\sigma)u \circ \phi_{\tau\sigma}^{-1} = 
^{\tau\sigma}u.$$
Here we need a little explanation for the second last equality:  
In fact, since $\bar{\phi}_{\tau} \circ \bar{\tau}\phi_{\sigma} = \bar{\phi}_{\tau\sigma}$, we can write 
$\phi_{\tau} \circ \tau\phi_{\sigma} = \phi_{\tau\sigma} \circ v$ by means of an element $v$ of 
$\mathrm{PAut}(\tau\sigma\eta ; id\vert_{\tau\sigma\bar{\eta}})$. 
Then we have 
$$\phi_{\tau} \circ \tau\phi_{\sigma} \circ (\tau\sigma)u \circ 
(\tau \phi_{\sigma})^{-1} \circ \phi^{-1}_{\tau} = \phi_{\tau\sigma} \circ v \circ 
(\tau\sigma)u \circ v^{-1} \circ \phi_{\tau\sigma}^{-1},$$ but since $\mathrm{PAut}(\tau\sigma\eta; id\vert_{\bar{\tau\sigma\eta}})$ is an 
abelian group, we see that $v \circ (\tau\sigma)u \circ v^{-1} = (\tau\sigma)u$. 

By the argument above, $\mathbf{C}^*$ acts on $\mathrm{PAut}(\eta; id\vert_{\bar{\eta}})$ from the left.

Next, for $\sigma, \tau \in \mathbf{C}^*$ we define $$f(\sigma, \tau) := \phi_{\sigma} \circ 
\sigma\phi_{\tau} \circ \phi_{\sigma\tau}^{-1} \in \mathrm{PAut}(\eta; id\vert_{\bar{\eta}}), $$
$$ \eta \stackrel{\phi_{\sigma\tau}^{-1}}\to \sigma\tau\eta \stackrel{\sigma\phi_{\tau}}\to 
\sigma\eta \stackrel{\phi_{\sigma}}\to \eta.$$
By using the fact that $\mathrm{PAut}(\eta; id\vert_{\bar{\eta}})$is abelian, we have, for  $\sigma, \tau, \rho \in 
\mathbf{C}^*$, 
$$f(\sigma\tau, \rho)\circ f(\sigma, \tau\rho)^{-1}\circ f(\sigma, \tau) = 
f(\sigma, \tau) \circ f(\sigma\tau, \rho)\circ f(\sigma, \tau\rho)^{-1}$$ 
$$ = (\phi_{\sigma}\circ \sigma\phi_{\tau}\circ \phi_{\sigma\tau}^{-1}) \circ 
(\phi_{\sigma\tau} \circ \sigma\tau\phi_{\rho} \circ \phi^{-1}_{\sigma\tau\rho}) \circ 
(\phi_{\sigma\tau\rho} \circ (\sigma \phi_{\tau\rho})^{-1} \circ \phi_{\sigma}^{-1})$$ 
$$ = \phi_{\sigma} \circ \sigma\phi_{\tau} \circ \sigma\tau\phi_{\rho} \circ 
\sigma\phi^{-1}_{\tau\rho} \circ \phi_{\sigma}^{-1} = \phi_{\sigma} \circ \sigma(\phi_{\tau} 
\circ \tau\phi_{\rho} \circ \phi_{\tau\rho}^{-1})  \circ \phi_{\sigma}^{-1}$$ 
$$ = {}^{\sigma}f(\tau, \rho).$$
This means that 
 $$f: \mathbf{C}^* \times \mathbf{C}^* \to  
\mathrm{PAut}(\eta; id\vert_{\bar{\eta}})$$ determines a 2-cocycle with respect to the group cohomology (Hochschild cohomology) for the rational representation $\mathrm{PAut}(\eta; id\vert_{\bar{\eta}})$ of the algebraic torus $\mathbf{C}^*$.  
Since an algebrac torus is linearly reductive, we have the vanishing of the higher Hochschild cohomology (cf. [Mi, Proposition 15.16]): 
$$H^i(\mathbf{C}^*, 
\mathrm{PAut}(\eta; id\vert_{\bar{\eta}})) = 0 \:\: i > 0.$$  
In particular, $f$ is a 2-coboundary. In other words, there exists a family of Poisson automorphisms 
$\{u_{\sigma}\}$, $u_{\sigma}  
\in \mathrm{PAut}(\eta; id\vert_{\bar{\eta}})$, parametrized by the elements of $\mathbf{C}^*$ such that the equalities holds: 
$$ f(\sigma, \tau) = {}^{\sigma}u_{\tau} \circ u_{\sigma\tau}^{-1} \circ u_{\sigma}.$$

Therefore, we have $$\phi_{\sigma}\circ \sigma\phi_{\tau} \circ \phi_{\sigma\tau}^{-1} = 
{}^{\sigma}u_{\tau} \circ u_{\sigma\tau}^{-1} \circ u_{\sigma}$$   
$$= {}^{\sigma}u_{\tau} \circ u_{\sigma} \circ u_{\sigma\tau}^{-1} = \phi_{\sigma} \circ 
\sigma u_{\tau} \circ  \phi_{\sigma}^{-1} \circ u_{\sigma} \circ u_{\sigma\tau}^{-1}.$$ 
Here, in the second equality,  we used the fact that $\mathrm{PAut}(\eta; id\vert_{\bar{\eta}})$ is abelian. 
 
Now, by operating $\phi_{\sigma}^{-1}$ from the left, and $\phi_{\sigma\tau}$ from the right, on both sides of this equality, we get 
$$\sigma \phi_{\tau} = \sigma u_{\tau} \circ \phi_{\sigma}^{-1} \circ u_{\sigma} \circ 
u_{\sigma\tau}^{-1} \circ \phi_{\sigma\tau}.$$ Further, by operating $\sigma(u_{\tau}^{-1}) = (\sigma u_{\tau})^{-1}$ from the left, 
on the both sides of this equality, we get  
$$\sigma(u_{\tau}^{-1}\circ \phi_{\tau}) 
= \phi_{\sigma}^{-1} \circ u_{\sigma} \circ u_{\sigma\tau}^{-1} \circ \phi_{\sigma\tau}.$$
Finally, operate $u_{\sigma}^{-1} \circ \phi_{\sigma}$ from the left on both sides. Then we have  
$$u_{\sigma}^{-1} \circ \phi_{\sigma} \circ \sigma(u_{\tau}^{-1}\circ \phi_{\tau}) =  
u_{\sigma\tau}^{-1} \circ \phi_{\sigma\tau}.$$ 
Therefore, if we put $\psi_{\sigma} := u_{\sigma}^{-1} \circ \phi_{\sigma}$ for each $\sigma \in \mathbf{C}^*$, then we have 
$$\psi_{\sigma}  \circ \sigma\psi_{\tau} = \psi_{\sigma\tau}.$$  
$\square$
\vspace{0.2cm}

The $\mathbf{C}^*$-action constructed in Theorem (2.1) is unique by the next proposition.  
\vspace{0.2cm}

{\bf Proposition (2.3)}. {\em Assume that, for a formal Poisson deformation $\{(X_n, \{\:, \:\}_n)\}$ of $(X, \{\:, \:\})$ (which is not  necessarily the universal formal Poisson deformation), there are two $\mathbf{C}^*$-actions satisfying the conditions (i), (ii) 
of Theorem (2.1), which induce the same $\mathbf{C}^*$-action on the base $\{S_n\}$. 
Then there exists an $\{S_n\}$-Poisson automorphism $h = \{h_n\}$ of $\{X_n\}$ such that $h_0 = id_X$ and 
it is equivariant with respect to the two $\mathbf{C}^*$-actions. }
\vspace{0.2cm}

{\bf Proof}. Suppose that a Poisson automorphism satisfying the conditions of the proposition is constructed up to $h_{n-1}$. 
By the prorepresentability of $\mathrm{PD}_X$, we can extend $h_{n-1}$ to a Poisson automorphism 
$h'_n$ of $X_n$. We compare the 1-st $\mathbf{C}^*$-action and the $\mathbf{C}^*$-action obtained by pulling back the 2-nd one 
by $h'_n$. As in the latter part of the proof of Theorem (2.1), we write $\eta$ for the pair of 
$\{X_0 \stackrel{i}\to X_n \stackrel{j}\to S_n\}$ and the Poisson structure $\{\:, \:\}_n$ on $X_n$. 

For $\sigma \in \mathbf{C}^*$, we consider the pair of $\{X_0 \stackrel{i\circ \sigma^{-1}}\to X_n \stackrel{\sigma \circ j}\to S_n\}$ and the Poisson structure $\sigma^{-l}\{\:, \:\}_n$ on $X_n$, which we denote by $\sigma\eta$. Note that, by the assumption of the proposition,  
$\sigma\eta$ is the same one, independent of the choice of  the two $\mathbf{C}^*$-actions. 

An isomorphism $\sigma\eta \stackrel{\phi_{\sigma}}\to \eta$ is determined by the 1-st $\mathbf{C}^*$-action. 
As in the proof of Theorem (2.1), we let $\mathbf{C}^*$ act on $\mathrm{PAut}(\eta ; \: id\vert\bar{\eta})$ by means of 
$\{\phi_{\sigma}\}$. For $u \in \mathrm{PAut}(\eta ; \: id\vert\bar{\eta})$, we denote by ${}^{\sigma}u$ the element obtained from 
$u$ by operating $\sigma$. 
 
On the other hand, an isomorphism $\sigma\eta \stackrel{\psi_{\sigma}}\to \eta$ is determined by the 2-nd $\mathbf{C}^*$-action.  
Then, by using an element $u_{\sigma}$ of $\mathrm{PAut}(\eta ; \: id\vert\bar{\eta})$, we can write $\psi_{\sigma} = u_{\sigma}^{-1} \circ \phi_{\sigma}$. Since 
$$\psi_{\sigma}  \circ \sigma\psi_{\tau} = \psi_{\sigma\tau},$$ we have 
$$u_{\sigma}^{-1} \circ \phi_{\sigma} \circ \sigma(u_{\tau}^{-1}\circ \phi_{\tau}) =  
u_{\sigma\tau}^{-1} \circ \phi_{\sigma\tau}.$$ 

From here, by tracing back the operation in the last part of the proof of Theorem (2.1), we get 
$$\phi_{\sigma}\circ \sigma\phi_{\tau} \circ \phi_{\sigma\tau}^{-1} = 
{}^{\sigma}u_{\tau} \circ u_{\sigma\tau}^{-1} \circ u_{\sigma}.$$
Since $\phi_{\sigma}\circ \sigma\phi_{\tau} \circ \phi_{\sigma\tau}^{-1} = 1$, 
$$u: \mathbf{C}^* \to \mathrm{PAut}(\eta ; \: id\vert\bar{\eta}), \:\:\: \sigma \to u_{\sigma}
$$ determines a 1-cocycle with respect to the group cohomology (Hochschild cohomology) for the rational representation 
$\mathrm{PAut}(\eta; id\vert_{\bar{\eta}})$ of the 
algebraic torus $\mathbf{C}^*$. Since an algebraic torus is linearly reductive,  we have $H^1(\mathbf{C}^*, 
\mathrm{PAut}(\eta; id\vert_{\bar{\eta}})) = 0$ (cf. [Mi, Proposition 15.16]).   
In particular, $\{u_{\sigma}\}$ is a 1-coboundary. In other words, we can write 
$u_{\sigma} = {}^{\sigma}\theta \circ \theta^{-1}$ for some $\theta \in \mathrm{PAut}(\eta; id\vert_{\bar{\eta}})$. 
Then we have $$\psi_{\sigma} = u_{\sigma}^{-1} \circ \phi_{\sigma} = \theta \circ {}^{\sigma}\theta^{-1} 
\circ \phi_{\sigma} = \theta \circ (\phi_{\sigma} \circ \sigma \theta^{-1} \circ \phi_{\sigma}^{-1}) \circ \phi_{\sigma} = \theta \circ \phi_{\sigma} \circ \sigma \theta^{-1}.$$
Namely, the following diagram commutes:   
\begin{equation} 
\begin{CD}
\sigma \eta @>{\phi_{\sigma}}>> \eta \\
@V{\sigma \theta}VV @V{\theta}VV \\ 
\sigma \eta @>{\psi_{\sigma}}>> \eta          
\end{CD} 
\end{equation}
Now, by putting $h_n := \theta \circ h'_n$, we have a Poisson automorphism $h_n$ which is a lifting of $h_{n-1}$ and satisfies the desired property of the proposition.  
$\square$ 

\vspace{0.2cm}

As a corollary of Proposition (2.3) we can prove that the universal formal Poisson deformation 
$\mathcal{X}^{univ} := \{(X_n, \{\:, \}_n)\}$ of $(X, \{\:, \:\})$ is also the universal $\mathbf{C}^*$-equivariant Poisson deformation of  
$(X, \{\:, \:\})$.

Let us make it more precise. 
Let $(A, m)$ be an Artinian local $\mathbf{C}$-algebra with 
$A/m = \mathbf{C}$. Moreover, we assume that $A$ has a $\mathbf{C}^*$-action. Let $(\mathrm{Art})_{\mathbf C}^{\mathbf{C}^*}$ 
be the category with these objects whose morphisms are $\mathbf{C}^*$-equivariant ones.    
We put $S := \mathrm{Spec}(A)$. 
A Poisson deformation $(\mathcal{X}, \{\:, \:\}_{\mathcal X}) \to S$ of $(X, \{\:, \:\})$ is called a $\mathbf{C}^*$-equivariant 
Poisson deformation if the underlying map $\mathcal{X} \to S$ is $\mathbf{C}^*$-equivariant and $\{\:, \:\}_{\mathcal X}$ has 
weight $-l$. Two $\mathbf{C}^*$-equivariant Poisson deformations $(\mathcal{X}, \{\:, \:\}_{\mathcal X}) \to S$ and 
$(\mathcal{X}', \{\:, \:\}_{{\mathcal X}'}) \to S$ are equivalent if there exsits a $\mathbf{C}^*$-equivariant 
Poisson $S$-isomorphism $(\mathcal{X}, \{\:, \:\}_{\mathcal X}) \cong (\mathcal{X}', \{\:, \:\}_{{\mathcal X}'})$ which restricts to 
the identity map of $X$. Let us consider the functor 
$$\mathrm{PD}_X^{\mathbf{C}^*}: (\mathrm{Art})_{\mathbf C}^{\mathbf{C}^*} \to (\mathrm{Set}) $$ 
which sends $A \in (\mathrm{Art})_{\mathbf C}^{\mathbf{C}^*}$ to the set of equivalence classes of $\mathbf{C}^*$-equivariant 
Poisson deformations of $(X, \{\:, \:\})$ over $S$. Now we have: 
\vspace{0.2cm}

{\bf Corollary (2.4)}. {\em The universal formal Poisson deformation $\mathcal{X}^{univ}$ of $(X, \{\:, \:\})$ is also the 
universal one among the $\mathbf{C}^*$-equivarint Poisson deformations.} 
\vspace{0.2cm}

{\bf Proof}. Let $R$ be the prorepresentable hull of $\mathrm{PD}_X$. 
If we are given a $\mathbf{C}^*$-equivariant Poisson deformation $(\mathcal{X}, \{\:, \:\}_{\mathcal X}) \to S$, then there 
ia a unique $\mathbf{C}^*$-equivariant map $\varphi: S \to \mathrm{Spec}(R)$ with $\varphi (0) = 0$. By definition, two Poisson deformations $\mathcal{X} \to S$ and $\mathcal{X}^{univ} \times_{\mathrm{Spec}(R)} S \to S$ are equivalent. In particular, there is a 
Poisson $S$-isomorphism $\Psi: \mathcal{X} \cong \mathcal{X}^{univ} \times_{\mathrm{Spec}(R)} S$. The right hand side admits a 
$\mathbf{C}^*$-action induced from the $\mathbf{C}^*$-action on $\mathcal{X}^{univ}$. We pull back this $\mathbf{C}^*$-action to 
a $\mathbf{C}^*$-action on $\mathcal{X}$ by $\Psi$. We compare this $\mathbf{C}^*$-action with the original $\mathbf{C}^*$-action 
on $\mathcal{X}$. But these two $\mathbf{C}^*$-actions are equivalent by Proposition (2.3). $\square$   
\vspace{0.2cm}

\S 3. {\bf Linearizations of line bundles}. 
\vspace{0.2cm}

The proof of Lemma A.8 of [Na] contains an error. 
Instead, we give here a different proof. More precisely, we prove Lemma (3.1), which is a slightly modified version of  [Na, Lemma A.8]. 
Lemma (3.1) is enough for the argument in [Na, Appendix] (see Remark below).

Assume that $(\hat{A}, m)$ is a complete local $\mathbf{C}$-algebra with $\hat{A}/m = \mathbf{C}$ and suppose that 
$\hat{A}$ has a $\mathbf{C}^*$-action. We put $\hat{Y} := \mathrm{Spec} (\hat{A})$. 
\vspace{0.2cm}

{\bf Lemma (3.1)}. 
{\em Let  $\hat{f}: \hat{X} \to \hat{Y}$ be a $\mathbf{C}^*$-equivarinat birational projective morphism. Assume
that $\hat{f}_*\mathcal{O}_{\hat X} = \mathcal{O}_{\hat Y}$. Let $\hat{L}$ be a line bundle on $\hat{X}$ such that 
$\hat{a}^*\hat{L} \cong pr_2^*\hat{L}$, where $\hat{a}: \mathbf{C}^* \hat{\times} \hat{X} \to \hat{X}$ is the $\mathbf{C}^*$-action and 
$pr_2: \mathbf{C}^* \hat{\times} \hat{X} \to \hat{X}$ is the projection. 
Then $\hat{L}$ has a $\mathbf{C}^*$-linearization.}
\vspace{0.2cm}

{\bf Proof}. {\bf (i)}  We put $Y_n := \mathrm{Spec} (A/m^{n+1})$ and $X_n := \hat{X} \times_{\hat Y}Y_n$. We denote by 
$L_n := \hat{L}\vert_{X_n}$.  
 
{\bf (ii)}  We restrict $\hat{a}$  to $\{1\} \times \hat{X}$. 
Then we get an automorphism of $\hat{L}$.  This automorphism can be written as 
$$\varphi_0 : \hat{L} \to \hat{L}, \:\: x \to \varphi_0 \cdot  x$$ by using an element $\varphi_0$ of 
$\Gamma(\hat{X}, \mathcal{O}_{\hat X})^*$.  
Consider here $pr_2^*\varphi_0^{-1}: pr_2^*\hat{L} \to pr_2^*\hat{L}$ and let us define an isomorphism 
$$\phi_0: \hat{a}^*\hat{L} \to 
pr_2^*\hat{L}$$ as a composition of the maps 
$\hat{a}^*\hat{L} \to pr_2^*\hat{L}$ and $pr_2^*\varphi_0^{-1}$. Then the restriction of $\phi_0$ to $1 \times \hat{X}$ is the identity map 
of $\hat{L}$. 
To measure how this isomorphism differes from a $\mathbf{C}^*$-linearization of 
$\hat{L}$, we prepare the following notation. 
For $\tau \in \mathbf{C}^*$, 
$\phi_0$ induces an isomorphism $\tau^*\hat{L} \to \hat{L}$. We denote it by $\phi_{0. \tau}$.
For $\sigma \in \mathbf{C}^*$, $\tau^*\hat{L} \stackrel{\phi_{0. \tau}}\to \hat{L}$ induces an isomorphism $\sigma^*\tau^*\hat{L} \to \sigma^*\hat{L}$. We denote it by $\sigma\phi_{0. \tau}$. 
If , for all $\sigma, \tau \in \mathbf{C}^*$, 
$$\phi_{0, \sigma} \circ \sigma \phi_{0, \tau} = \phi_{0, \sigma\tau}$$
hold, then $\phi_0$ gives a $\mathbf{C}^*$-linearization of  $\hat{L}$. 
We put here $$f(\sigma, \tau) := \phi_{0, \sigma} \circ \sigma \phi_{0, \tau} \circ \phi_{0, \sigma\tau}^{-1}.$$
Then $f(\sigma, \tau)$ gives an automorphism of $\hat{L}$. An automorphism of $\hat{L}$ determines an element 
of $\Gamma (\hat{X}, \mathcal{O}_{\hat X})^*$. 
Then an element $F$ of $\Gamma ((\mathbf{C}^*)^2 \hat{\times} \hat{X}, \mathcal{O}_{(\mathbf{C}^*)^2 \hat{\times} \hat X})$ is determined and if we substitute $F$ a particular point $(\sigma, \tau)$ of $\mathbf{C}^*$, then we get 
$f(\sigma, \tau)$. Since one can write  
$$\Gamma ((\mathbf{C}^*)^2 \hat{\times} \hat{X}, \mathcal{O}_{(\mathbf{C}^*)^2 \hat{\times} \hat{X}}) = \Gamma ((\mathbf{C}^*)^2 \hat{\times} \hat{Y}, \mathcal{O}_{(\mathbf{C}^*)^2 \hat{\times} \hat X})$$
$$= \mathbf{C}[s, t, 1/s, 1/t] \hat{\otimes} \hat{A} = \lim_{\leftarrow} (\mathbf{C}[s, t, 1/s,  1/t]\otimes A/m^{k+1}), $$ we have $$F = \lim_{\leftarrow} f_k,\:\:  f_k \in \mathbf{C}[s, t, 1/s,  1/t]\otimes A/m^{k+1}$$ and $f_0 \in \mathbf{C}[s, t, 1/s, 1/t]^*$. 
Hence we can write as $f_0 = c s^at^b$, $c \in \mathbf{C}^*$, $a, b  \in \mathbf{Z}$. 
By the construction of $\phi_0$, we get $f(\sigma, 1) = f(1, \tau) = 1$ for arbitrary $\sigma, \tau \in \mathbf{C}^*$. This implies that $f_0 = 1$. 

{\bf (iii)} 
The rough idea of the proof is the following. First choose an element $u_1 \in 1 + \mathbf{C}[t, 1/t] \hat{\otimes} m\hat{A}$ so that the isomorphism $\phi_1\vert_{{\hat a}^*L_1} : {\hat a}^*L_1 \to pr_2^*L_1$ defined as the composite of 
$u_1: pr_2^*\hat{L} \to pr_2^*\hat{L}$ and $\phi_0$ gives a $\mathbf{C}^*$-linearization of $L_1$. Next choose $u_2 \in 1 + \mathbf{C}[t, 1/t] \hat{\otimes} m^2\hat{A}$ so that  $\phi_2 := u_2 \circ \phi_1$ gives a $\mathbf{C}^*$-linearization of $L_2$. Repeating this procedure, we finally put  $u_{\infty} := \prod u_k$ and consider 
$\phi := u_{\infty} \circ \phi_0$. Then $\phi$ gives a $\mathbf{C}^*$-linearization of 
$\hat{L}$.  

We prove it by the induction on $k$.  Suppose now that we can choose $u_{k-1}$ so that $\phi_{k-1}$ gives a $\mathbf{C}^*$-linearization of $L_{k-1}$ (when $k = 1$, we put $u_0 = 1$).  
The homomorphism $\phi_{k-1}: \hat{a}^*\hat{L} \to pr_2^*\hat{L}$ induces, for $L_k$ and  $L_{k-1}$,  
two homomorphisms $$\bar{\phi}'_k : \hat{a}^*L_k \to pr_2^*L_k, \:\:\: \bar{\phi}_{k-1}: 
{\hat a}^*L_{k-1} \to pr_2^*L_{k-1}.$$ Moreover, for 
$\sigma \in \mathbf{C}^*$, $\bar{\phi}'_k$ and $\bar{\phi}_{k-1}$ determine homomorphisms  
$$\bar{\phi}'_{k, \sigma}: \sigma^*L_k \to L_k, \:\: \bar{\phi}_{k-1, \sigma} : \sigma^*L_{k-1} 
\to L_{k-1}.$$ Here we put $$f(\sigma, \tau) := {\bar{\phi}'}_{k, \sigma} \circ \sigma {\bar{\phi}'}_{\tau} \circ 
\bar{\phi'}^{-1}_{k, \sigma\tau}$$ for $\sigma, \tau \in \mathbf{C}^*$. 
Then this is an automorphism of $L_k$. Since 
$\bar{\phi}_{k-1}$ gives a $\mathbf{C}^*$-linearization of  $L_{k-1}$, this automorphism restricts to the identity map  
of $L_{k-1}$. Hence we have $f(\sigma, \tau) \in \mathrm{Aut}(L_k ; id\vert_{L_{k-1}})$. Here let us consider the 
$\mathbf{C}$-linear subspace $$id + m^k/m^{k+1} \subset \mathrm{Aut}(L_k ; id\vert_{L_{k-1}})$$ of 
$\mathrm{Aut}(L_k ; id\vert_{L_{k-1}})$. Noticing that $\bar{\phi}'_k$ is obtained from $\phi_{k-1}$, we see that 
$$ f(\sigma, \tau) \in id + m^k/m^{k+1}.$$  
Moreover, for $v \in id + m^k/m^{k+1}$, we define 
$${}^{\sigma}v := \bar{\phi'}_{k, \sigma} \circ \sigma v \circ \bar{\phi'}^{-1}_{k, \sigma}.$$
Then we have ${}^{\sigma}v  \in id + m^k/m^{k+1}$, by which 
$\mathbf{C}^*$ acts on $id + m^k/m^{k+1}$. 
Now $$f : \mathbf{C}^* \times \mathbf{C}^* \to id + m^k/m^{k+1}$$ is a 2-cocyle for the Hochschild cohomology for 
the rational representation $id + m^k/m^{k+1}$ of the algebraic torus $\mathbf{C}^*$.      
Since higher Hochschild cohomology vanishes for an algebraic torus (cf. [Mi, Proposition 15.16]), $f$ is a 2-coboundary.     
In other words, there is a family of elements $\{\bar{v}_{k, \sigma}\}$, $\bar{v}_{k,   \sigma} \in id + m^k/m^{k+1} $ 
parametrized by elements $\sigma$ of 
$\mathbf{C}^*$ so that  
$$ f(\sigma, \tau) = {}^{\sigma}\bar{v}_{k, \tau} \circ \bar{v}_{k, \sigma\tau}^{-1} \circ \bar{v}_{k, \sigma}. $$

If we reset $\bar{\phi}_{k, \sigma}$ by 
$\bar{\phi}_{k, \sigma} := \bar{v}_{k, \sigma}^{-1} \circ \bar{\phi}'_{k, \sigma}$, then we have 
$$\bar{\phi}_{k, \sigma} \circ \sigma \bar{\phi}_{k, \tau} = \bar{\phi}_{k, \sigma\tau }.$$ 
Since $\{\bar{v}_{k. \sigma}\}$ depends algebraically on $\sigma$, an element 
$$ \bar{v}_k \in 1 + \mathbf{C}[t, 1/t] \otimes m^k/m^{k+1}$$
is determined. If we substitute $t = \sigma$, then we get $\bar{v}_{k, \sigma}$. 
Here put $\bar{u}_k := \bar{v}_k^{-1} \in 1 + \mathbf{C}[t, 1/t] \otimes m^k/m^{k+1}$ and take a lift 
$u_k \in 1 +  \mathbf{C}[t, 1/t] \hat{\otimes} m^k\hat{A}$
of $\bar{u}_k$. Then this $u_k$ is the desired element.  
$\square$ \vspace{0.2cm}

{\bf Remark}.  In [Na], Lemma (3.1) is used in the following context. 
Let $Y := \mathrm{Spec}(A)$ be an affine normal variety with rational singularities.  
We further assume that $Y$ has a good $\mathbf{C}^*$-action with a unique fixed point $0 \in X$.
Here let us asuume that $f: X \to Y$ is a $\mathbf{C}^*$-equivaraint partial resolution of $Y$.  In particular, 
$X$ is normal. 
Let $m$ be the maximal ideal of $A$ corresponding to $0$, and let $\hat{A}$ be the completion of $A$ along $m$.
Put $\hat{Y} := \mathrm{Spec}(\hat{A})$ and $\hat{X} := X \times_X \hat{Y}$. We denote by $\hat{f}: \hat{X} \to \hat{Y}$ 
the natural projective morphism induced by $f$. Then this $\hat{f}$ is nothing but the $\hat{f}$ in Lemma (3.1). 
Let $\mathcal{L}$ be a line bundle on $X^{an}$. Put 
$\mathcal{L}_n := \mathcal{L}\vert_{X_n^{an}}$. Since $X_n$ is proper over $\mathbf{C}$, there exists a unique algebraic line bundle 
$L_n \in \mathrm{Pic}(X_n)$ such that 
$L_n^{an} = \mathcal{L}_n$ by the GAGA principle. 
By the Grothendieck existence theorem ([EGA III] 1, 5.1.6), we have 
$$\mathrm{Pic}(\hat{X}) \cong \lim_{\leftarrow}\mathrm{Pic}(X_n).$$ 
Hence $\{L_n\}$ determines a line bundle $\hat{L}$ on $\hat{X}$. This is nothing but the line bundle $\hat{L}$ in 
Lemma (3.1). 

We write the $\mathbf{C}^*$-action on $X$ as 
$a: \mathbf{C}^* \times X \to X$. 
$a$ determines a $\mathbf{C}^*$-action 
$\hat{a}: \mathbf{C}^* \hat{\times} 
\hat{X} \to \hat{X}$ on $\hat{X}$.  Let us prove \vspace{0.2cm}

{\bf Claim}. {\em $\hat{a}^*\hat{L} \cong pr_2^*\hat{L}$.
In other words, the assumption of Lemma (3.1) holds.}  \vspace{0.2cm}

{\em Proof}.  
$a$ induces a $\mathbf{C}^*$-action on $X^{an}$: 
$a^{an}: \mathbf{C}^* \times X^{an} \to X^{an}$. 
Since $H^i(\mathbf{C}^* \times X^{an}, \mathcal{O}_{\mathbf{C}^* \times X^{an}}) = 0$ $(i = 1, 2)$, 
we have $\mathrm{Pic}(\mathbf{C}^* \times X^{an}) \cong H^2(\mathbf{C}^* \times X^{an}, \mathbf{Z})$. 
Since $H_*(\mathbf{C}^*, \mathbf{Z})$ is a free module, the Kunneth formula yields $$H^2(\mathbf{C}^* \times X^{an}, \mathbf{Z}) = H^0(\mathbf{C}^*, \mathbf{Z})\otimes 
H^2(X^{an}, \mathbf{Z}) \oplus H^1(\mathbf{C}^*, \mathbf{Z})\otimes H^1(X^{an}, \mathbf{Z}).$$
We notice here that $H^1(X^{an}, \mathbf{Z}) = 0$. In fact, since $Y^{an}$ has only rational singularities, 
one has $R^1f^{an}_*\mathbf{Z} = 0$. Then the Leray spectral sequence yields the exact sequence 
$$ 0 \to H^1(Y^{an}, \mathbf{Z}) \to H^1(X^{an}, \mathbf{Z}) \to 
H^0(Y^{an}, R^1f^{an}_*\mathbf{Z}).$$ Since $Y^{an}$ is conical, it can be contractible to the origin and, hence we get 
$H^1(Y^{an}, \mathbf{Z}) = 0$. 
Therefore $$pr_2^*: \mathrm{Pic}(X^{an}) \to \mathrm{Pic}(\mathbf{C}^* \times X^{an})$$ is an 
isomorphism. This means that $$(a^{an})^*\mathcal{L} \cong pr_2^*\mathcal{L}.$$ In fact, 
since $pr_2^*$ is an isomorphism, one can write 
$(a^{an})^*\mathcal{L} = pr_2^*K$, $K \in \mathrm{Pic}(X^{an})$. 
As  $(a^{an})^*\mathcal{L}\vert_{1 \times X^{an}} = \mathcal{L}$, we see that 
$K = \mathcal{L}$.     
The $\mathbf{C}^*$-action $a^{an}$ induces a $\mathbf{C}^*$-action on $X_n^{an}$ for each $n \geq 0$. 
We denote it by the same $a^{an}$.  We put $\mathcal{L}_n := \mathcal{L}\vert_{X_n}$.
Then the isomorphism above induces an isomorphism $(a^{an})^*\mathcal{L}_n \cong pr_2^*\mathcal{L}_n$ of line bundles on 
$\mathbf{C}^* \times X_n^{an}$. As we remarked at the beginning,  
$L_n \in \mathrm{Pic}(X_n)$ is an algebraic line bundle such that  $L_n^{an} = \mathcal{L}_n$.  
Let us prove that $a^*L_n \cong pr_2^*L_n$. 
In order to do so, we put $W_n := \mathrm{Spec}\: \Gamma (X_n, \mathcal{O}_{X_n})$ and denote by $g_n$ by the natural morphism 
$X_n \to W_n$. Then $id \times f_n$ factorizes as follows:  
$$ \mathbf{C}^* \times X_n \stackrel{id \times g_n}\to \mathbf{C}^* \times W_n 
\to \mathbf{C}^* \times Y_n.$$
Since $(a^{an})^*L_n^{an} \otimes (pr_2^*L_n^{an})^{-1}$ is trivial, we have  
$$\mathcal{O}_{\mathbf{C}^* \times W_n^{an}} = (id \times g_n)^{an}_*
[(a^{an})^*L_n^{an} \otimes (pr_2^*L_n^{an})^{-1}]$$ $$= [(id \times g_n)_*(a^*L_n 
\otimes (pr_2^*L_n)^{-1})]^{an}.$$ 
It follows from this fact that $M := (id \times g_n)_*(a^*L_n 
\otimes (pr_2^*L_n)^{-1})$ is a line bundle on $\mathbf{C}^* \times W_n$. 
On the other hand, since $W_n$ is a local Artinian $\mathbf{C}$-scheme and $\mathrm{Pic}(\mathbf{C}^*) 
= 1$, we can show that  $\mathrm{Pic}(\mathbf{C}^* \times W_n) = 1$. 
In order to do that, viewing $\mathbf{C}^* \times W_n$ as a deformation of $\mathbf{C}^*$ over $W_n$, we represent   
$W_n$ as a sequence of small cloed immersion $$S_0 := \mathrm{Spec} \: \mathbf{C} \to S_1 \to S_2 
\to ... \to S_m = W_n.$$ This determines a sequence of closed immersion 
$$\mathbf{C}^* \times S_0 \to \mathbf{C}^* \times S_1 \to \mathbf{C}^* \times S_2 \to ... 
\to \mathrm{C}^* \times S_m.$$
Since  $H^i(\mathbf{C}^*, \mathcal{O}_{\mathbf{C}^*})
= 0$, $i = 1, 2$, the restriction maps of Picard groups
$$\mathrm{Pic}(\mathbf{C}^* \times S_m) \to ... \to \mathrm{Pic}(\mathbf{C}^* \times 
S_2) \to \mathrm{Pic}(\mathbf{C}^* \times S_1) \to \mathrm{Pic}(\mathbf{C}^* \times S_0)$$ are all isomorphisms. Hence $M$ is a trivial line bundle. Moreover, since 
$a^*L_n \otimes (pr_2^*L_n)^{-1} = (id \times g_n)^*M$, we see that $a^*L_n \otimes (pr_2^*L_n)^{-1}$
is also trivial.  Again, by the Grothendieck existence theorem, we have  
$$\mathrm{Pic}(\mathbf{C}^* \hat{\times} \hat{X}) \cong 
\lim_{\leftarrow} \mathrm{Pic}(\mathbf{C}^* \times X_n).$$
As a consequence, we have shown that $\hat{a}^*\hat{L} \cong pr_2^*\hat{L}$. 
$\square$ 
\vspace{0.2cm}

\begin{center}
{\bf References} 
\end{center}
\vspace{0.2cm}

[EGA III] Grothendieck, A., Dieudonne, J.: 
Elements de Geometrie Algebrique, 
Etude cohomologique des faisceaux coherents, Publ. Math. IHES {\bf 11} (1961), {\bf 17}(1963)

[KKMS] Kempf, G.,  Knudsen, F., Mumford, D., Saint-Donat, B.: Toroidal embeddings. I. Lecture Notes in Mathematics, {\bf 339}. Springer-Verlag, Berlin-New York, 1973. viii+209 pp. 

[Mi] Milne, J.: Algebraic groups, Cambridge Studies in Advanced Mathematics, {\bf 170} Cambridge University Press, Cambridge, 2017. xiv+644 pp.

[Na] Namikawa, Y.: Flops and Poisson deformations of symplectic varieties, 
Publ. Res. Inst. Math. Sci. {\bf 44} (2008), 259-314

[Na 2] Namikawa, Y.: Poisson deformations of affine symplectic varieties, Duke Math. J. {\bf 156}, No.1 (2011), 51 - 85 

[R] Rim, D.: Equivariant $G$-structure on versal deformations, Trans. Amer. Math. Soc. {\bf 257} (1980), 217 - 226
\vspace{0.2cm}

\begin{center}
Research Institute for Mathematical Sciences, Kyoto University, Oiwake-cho, Kyoto, Japan

E-mail address: namikawa@kurims.kyoto-u.ac.jp  
\end{center}

\end{document}